\documentclass[%
 reprint,
nobibnotes,
 amsmath,amssymb,
 aps,
 prl,
]{revtex4-2}

\usepackage{graphicx}
\usepackage{dcolumn}
\usepackage{bm}
\usepackage{enumitem}
\usepackage{amssymb}
\usepackage{amsmath}
\usepackage{psfrag}
\usepackage{mathtools}
\usepackage{multirow}
\usepackage{color}
\usepackage{hyperref}
\usepackage[export]{adjustbox}
\usepackage{url}

\definecolor{cyan1}{rgb}{0.387, 0.82, 1}
\definecolor{red1}{rgb}{0.902, 0.383, 0.355}
\definecolor{RGBred}{rgb}{1,0,0}
\definecolor{RGBblue}{rgb}{0,0,1}
\definecolor{outgreen}{rgb}{0.38,0.73,0.66}
\definecolor{wavered}{rgb}{0.933,0.196,0.18}
\definecolor{waveblue}{rgb}{0.133,0.337,0.651}
\definecolor{wallgrey}{rgb}{0.4,0.4,0.4}
\definecolor{weakblue}{rgb}{0.757,0.792,0.89}
\definecolor{weakred}{rgb}{0.973,0.745,0.745}
\definecolor{red2}{cmyk}{0,0.2,0,0}

\newcommand{\tr}[1]{{#1}} 
\newcommand{\sr}[1]{{#1}} 
\newcommand{\rr}[1]{{#1}} 
\newcommand{\pc}[1]{{#1}} 
\newcommand{\fc}[1]{{#1}} 

\begin{document}

\title{Superradiant Scattering from Nonlinear \pc{Wave-Mode} Coupling}%

\author{Tiemo Pedergnana}
 \email{ptiemo@ethz.ch}
 \affiliation{%
 CAPS Laboratory, Department of Mechanical and Process Engineering, ETH Z{\"u}rich, Sonneggstrasse 3,
8092 Z{\"u}rich, Switzerland
}%
\author{Nicolas Noiray}%
 \email{noirayn@ethz.ch}
\affiliation{%
 CAPS Laboratory, Department of Mechanical and Process Engineering, ETH Z{\"u}rich, Sonneggstrasse 3,
8092 Z{\"u}rich, Switzerland
}%


\begin{abstract}
\pc{Waves scattered at a self-oscillating mode can exhibit superradiance, or net ampliﬁcation of an external harmonic excitation. This exotic behavior, arising from the nonlinear coupling between the mode and the incident wave, is theoretically predicted and experimentally conﬁrmed \fc{in this work}. We propose a generic theory of nonlinear wave-mode coupling, which is derived in analogy to the temporal coupled-mode theory of [Fan et al., J. Opt. Soc. Am. A 20, 569 (2003)]. A well-reproducible aeroacoustic realization of a superradiant scatterer was used to test the theory’s predictions. It is shown that the nonlinear wave-mode coupling can be exploited to quasi-passively tune the reflection and transmission coefficients of a side cavity in a waveguide. The theoretical framework used to describe this type of superradiance is applicable to non-acoustic systems and may be used to design lossless scattering devices.}
\end{abstract}
\date{\today}%
\maketitle

Superradiance occurs when an energetic, radiating region \tr{scatters incident waves with} a net amplification \cite{Vicente2018}. Paradigmatic manifestations are the Penrose process, by which particles gain momentum from a rotating black hole \cite{penrose1971extraction} and the Zel'dovich mechanism, where a rotating metallic cylinder increases electromagnetic oscillations \cite{zel1972amplification}. Extraction of energy from rotating systems has been \tr{investigated} in optics \cite{Belgiorno2010,Drori2019}, acoustics \cite{Cromb20201069}, hydrodynamics \cite{Acheson1976433,Torres2017833} and quantum mechanics \cite{Steinhauer2014864,Steinhauer2016959,MunozdeNova2019688}.

\begin{figure}[t!]
\begin{psfrags}
    
\psfrag{c}{\hspace{-0.03cm}$a(t)$\hspace{3.27cm}$a(t)$}
\psfrag{a}{\hspace{-0.35cm}\textcolor{waveblue}{$|s_\mathrm{in}\rangle_j$}\hspace{3.08cm}\textcolor{weakblue}{$|s_\mathrm{in}\rangle_j$}}
\psfrag{b}{\hspace{-0.35cm}\textcolor{weakred}{$|s_\mathrm{out}\rangle_j$}\hspace{2.945cm}\textcolor{wavered}{$|s_\mathrm{out}\rangle_j$}}
\psfrag{d}{\hspace{-0.24cm}NLC off \hspace{2.45cm} NLC on}
\psfrag{e}{\textcolor{wallgrey}{\hspace{-0.8cm}$j$th port}}

    \includegraphics[width=0.42\textwidth]{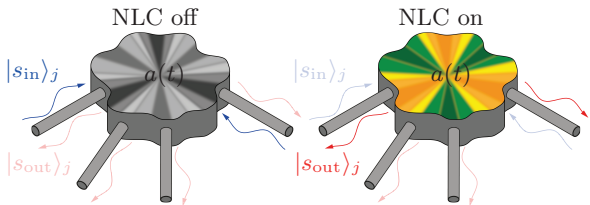}
    \end{psfrags}
            \vspace{-0.3cm}
    \caption{\pc{In a classic resonant cavity with multiple ports (left inset), due to dissipative effects, the energy of the incident waves $\langle s_\mathrm{in}|s_\mathrm{in}\rangle$ always exceeds that of the outgoing waves $\langle s_\mathrm{out}|s_\mathrm{out}\rangle$. This limitation is overcome by inducing synchronization between the natural radiation of a self-oscillating cavity mode and the incident waves, which can lead to superradiance, or net amplification of an external harmonic excitation (right inset).} The \pc{incident and outgoing waves in} the $j$th port are denoted by $|s_\mathrm{in}\rangle_j$ and $|s_\mathrm{out}\rangle_j$.}
    \label{Figure 1}
\end{figure}

\tr{Here,} we consider a scenario where \tr{superradiant scattering is achieved by using a self-oscillating cavity mode. More specifically,} energy is \tr{steadily} supplied to the \tr{cavity}, destabilizing \tr{one of its modes to small perturbations, thus leading, under the action of internal nonlinearities, to a} limit cycle \rr{$a_0=|a_0|\, e^{i\omega_0 t}$} at the self-sustained \rr{angular} frequency \tr{$\omega_0$.} If an external wave with frequency \tr{$\omega\approx \omega_0$} impinges on \pc{the cavity}, nonlinear \pc{wave-mode} coupling (NLC) takes place\pc{, leading to an amplitude-dependent} scattering matrix $S$ \pc{exhibiting superradiance (Fig. \ref{Figure 1}). The entries of $S$ may be tuned by varying the forcing amplitude}. 

\sr{To} model this phenomenon\sr{, we derive a scattering theory which is related to} the temporal coupled-mode theory (TCMT) of Fan et al. \cite{Fan2003569} describing unitary multiport \pc{scattering by} linearly stable resonant cavities \cite{Suh20041511,Wang2005369,Dmitriev20181165,Zhao2019}. In contrast, we are interested in modeling \pc{\textit{nonlinear}}, non-unitary \tr{scattering by} cavities \tr{exhibiting a self-oscillating} mode. Despite fruitful recent applications of coupled-mode theories in optics \cite{Sweeney2019,Li2020}, electrodynamics \cite{Alpeggiani2017,Mazor2021}, quantum mechanics \cite{Zhang2019,Franke2021,Ren2021} and acoustics \cite{Fleury2014516,Achilleos2017,Pagneux22}, no analogous framework exists for nonlinear systems. \pc{This work aims to fill this gap in the literature. For simplicity, we restrict our discussion to the case of a single self-oscillating mode in a cavity with two ports.}

The following set of equations is then used to model the nonlinear \sr{scattering process}: 
\begin{eqnarray}
    \dot{a}&=&(i\omega_0 +\nu -\kappa |a|^2 )a +D^\dag \tr{|s_\mathrm{in}\rangle},\label{modal dynamics}\\
   |s_\mathrm{out}\rangle&=&\underbrace{(C+D F^{-1} D^\dag)|s_\mathrm{in}\rangle}_{\pc{\text{linear}}} + \underbrace{D a}_{\pc{\text{nonlinear}}}, \label{input output relation}
\end{eqnarray}
\noindent\pc{where $a\in\mathbb{C}$ is the complex modal amplitude, $\dot{()}$ is the derivative with respect to time $t$, $\omega_0\in\mathbb{R}^+$ is the eigenfrequency, $\nu\in\mathbb{R}$ is the linear growth or decay rate, $\kappa\in\mathbb{R}^+$ is the nonlinearity constant, $D\in\mathbb{R}^2$ is the coupling matrix, $(\cdot)^\dag$ is the Hermitian conjugate, $\vert s_\mathrm{in}\rangle$ and $\vert s_\mathrm{out}\rangle\in\mathbb{C}^2$ are the incident and outgoing waves at the forcing frequency $\omega\in\mathbb{R}^+$, $C+D F^{-1} D^\dag\in\mathbb{C}^{2\times 2}$ is the frequency-dependent background scattering matrix, $C\in\mathbb{R}^{2\times 2}$ is the constant component of the background and $F\in\mathbb{C}$ is defined as $F=i(\omega-\omega_0)+\gamma$. Self-oscillation can occur only for $\nu>0$, which is the case considered here.}

\pc{The parameter $\gamma\in\mathbb{R}^+$ is the decay rate characteristic of the background process for a linearly stable resonant cavity in the absence of self-oscillation and internal nonlinearities:
\begin{eqnarray}
     \dot{a}&=&(i\omega_0-\gamma)a+D^\dag \vert s_\mathrm{in}\rangle, \label{linearized modal dynamics}\\
     \vert s_{\mathrm{out}}\rangle&=&C \vert s_{\mathrm{in}}\rangle+D a. \label{linearized input-output}
\end{eqnarray}
Equations \eqref{linearized modal dynamics} and \eqref{linearized input-output} are considered here to explain the background matrix $C+DF^{-1}D^\dag$. Together, they form the classic TCMT equations. Analogous to the constant matrix $C$ in TCMT models, the frequency-dependent matrix $C+DF^{-1}D^\dag$ represents the background of the nonlinear process in our theory described by Eqs. \eqref{modal dynamics} and \eqref{input output relation}. This matrix, which is derived in the Supplemental Material \cite{suppmat}, ensures linear resonant scattering for large forcing amplitudes, i.e., it bounds the nonlinear scattering process to incident waves of relatively small to moderate amplitudes \footnote{Note that in the steady state, due to the cubic nonlinearity in Eq. \eqref{modal dynamics}, the linear term dominates the nonlinear contribution to the outgoing wave \eqref{input output relation} for large forcing amplitudes \pc{$\Vert C\vert s_\mathrm{in}\rangle \Vert\gg \Vert D a_{0}\Vert$}.}. A physical explanation for such a model in the specific example examined in this study is given below.}

\pc{Here, we explore the theoretical framework defined by Eqs. \eqref{modal dynamics} and \eqref{input output relation} for the special case of an asymmetric, reflectional superradiant scatterer coupled to a one-dimensional (1D) waveguide. This choice is motivated by the nature of the present experimental realization of superradiant scattering, which is based on a classic type of aeroacoustic instability \cite{Howe1980407,Elder1982532} that has recently been revisited both theoretically and experimentally \cite{Boujo2020,Bourquard2020,Pedergnana2021}, and is accordingly well understood. Simply put, when a low-Mach air flow is incident on a cavity's aperture at the right angle and velocity, it can destabilize one of its resonant modes, which leads to self-sustained oscillations. This transition from resonance to whistling, which can be modeled as the supercritical Hopf bifurcation associated with a Stuart--Landau oscillator \cite{Stuart19581,landau1959fluid,arnold2012geometrical,Bourquard2020}, is exploited in this work to create a limit cycle in the cavity, which enables superradiant scattering of harmonic waves.} 

\pc{In Eq. \eqref{input output relation}, the resonant nature of the background matrix $C+DF^{-1}D^\dag$ governing scattering at large forcing amplitudes is consistent with experiments \cite{Bourquard2020} and numerical simulations \cite{Boujo2018386} on configurations similar to that of the experiments presented below: In \cite[Fig. 8]{Bourquard2020}, measurements showed that the resistance of the cavity-waveguide interface in the presence of mean flow, which can take negative values when the aeroacoustic feedback is constructive, exhibits a nonlinear behavior with respect to the forcing amplitude as its value tends towards the positive constant resistance measured in the absence of flow. This saturated and \textit{finite} positive resistance explains the resonant nature of the scattering dominated by the linear term in Eq. \eqref{input output relation} when the forcing amplitude is sufficiently large.}

\pc{In this work, the modeled linear background contains an asymmetry, the bias induced by a non-negligible mean flow inside the waveguide. Let us now discuss how we construct the associated matrices. Considering that previous works on TCMT constrained the links between $\gamma$, $C$ and $D$ on the basis of symmetry arguments \cite{Fan2003569,Wang2005369,Zhao2019} which are inherently difficult to generalize to non-symmetric systems, we developed another approach for deducing the coupling matrix $D$ from two user-defined quantities: The target matrix $S_*$ corresponding to optimal, resonant scattering and the matrix $C$, which is the frequency-independent component of the background scattering matrix. This construction of the coupling matrix can thus be used to customize the cavity-waveguide interface. In the superradiant model given by Eqs. \eqref{modal dynamics} and \eqref{input output relation}, the resulting coupling matrix $D$ also defines the contribution of the nonlinear mode $a$ to the outgoing wave \eqref{input output relation}.}

\pc{It is worth noting that the present modeling framework is by no means restricted to the example of 1D reflectional scattering considered here and may be extended to other systems in a straightforward manner by changing $S_*$ and $C$, as exemplified in the Supplemental Material \cite{suppmat}.}

We follow the convention by which the diagonal and off-diagonal elements of $S$ are reflection and transmission coefficients, respectively. To theoretically model the linear part of the scattering process, we superimpose a diagonal (reflective) target matrix $S_*$ and a purely transmissive constant background $C$:
\begin{eqnarray}
    S_*= \begin{pmatrix}
1+\epsilon/\sigma & 0\\
0 & 1-\epsilon/\sigma
\end{pmatrix},\quad
   C= \begin{pmatrix}
0 & 1\\
1 & 0
\end{pmatrix},
\end{eqnarray}
\pc{where $\epsilon$ is the asymmetry and $\sigma\in \mathbb{R}^+$ is the unitarity factor, which we introduced with the objective of accounting for internal losses, as is explained below. Formally approximating the scaled target matrix $\sigma S_*$ by the background matrix $C+DF^{-1}D^\dag$ yields a low-rank approximation problem from which the unknown coupling matrix $D$ can be determined:}
\begin{eqnarray}
     \sigma S_*-C\approx D F^{-1} D^\dag. \label{Low-rank approximation}
\end{eqnarray}
\pc{The eigenvalues $\mu_j$ and corresponding normalized eigenvectors $v_j$, $j=1,2$, of the matrix $\sigma S_*-C$ are given in the Supplemental Material \cite{suppmat}.}
 \pc{To capture the main characteristics of the specified matrices $\sigma S_*$ and $C$, we define the coupling matrix $D$ such that for resonant forcing ($F=\gamma$), the singular matrix $DF^{-1}D^\dag$ shares with $\sigma S_*-C$ its dominant subspectrum, i.e., the eigenvalue with the largest absolute value and the corresponding eigenvector. Both $D$ as well as the ideal parameter configuration, which we call ``perfect matching'', are defined by the following two consistency conditions:} (I) The \pc{spectrum of the singular matrix} $D\gamma^{-1}D^\dag$ must coincide with the \pc{dominant subspectrum} of $\sigma S_*-C$. (II) For a perfectly matched scatterer \pc{with} $\sigma=1$, $C+D\gamma^{-1}D^\dag=S_*$ \pc{must be} unitary.
 
 \rr{By the Hermiticity of the chosen $S_*$ and $C$, we can expand the $2$-by-$2$ matrix $\sigma S_*-C$ in terms of its eigenvectors $v_j\in\mathbb{C}^2$ and eigenvalues $\mu_j\in\mathbb{R}$ as follows:
\begin{eqnarray}
\sigma S_*-C=\mu_1 v_1 v_1^\dag + \mu_2 v_2 v_2^\dag,
\end{eqnarray}
where $|\mu_1|\geq|\mu_2|$. To satisfy (I), we define $D=v_1 \sqrt{\gamma \mu_1}$, so that $D \gamma^{-1} D^\dag=\mu_1 v_1 v_1^\dag$.} From (II), we deduce that the perfect matching condition is $\epsilon=0$. Written out, the coupling matrix $D$ is given by
\begin{eqnarray}
    D= \sqrt{\gamma h(\sigma,\epsilon)}  \begin{pmatrix}
-g(\epsilon) \\
1 
\end{pmatrix}, \label{coupling matrix D}
\end{eqnarray}
where $g(\epsilon)=\epsilon + \sqrt{\epsilon^2 + 1}$ and $h(\sigma,\epsilon)=(\sigma + \sqrt{\epsilon^2+1})/\big[g(\epsilon)^2+1\big]$. \pc{As we show in the Supplemental Material \cite{suppmat}, the asymmetry $\epsilon$ induces a biased gain in the system, which is also observed in the experiments presented below due to the presence of a mean flow with low, but non-negligible Mach number in the waveguide. Consequently, both the linear and the nonlinear part of the scattering matrix contribute to the reflectional superradiance.}

\pc{Following \cite{Zhao2019}, we account for dissipative losses in the linear background process with the internal decay rate $\gamma_\mathrm{i}$, which is related to the coupling matrix $D$: $\sr{D^\dag D = 2(\gamma-\gamma_\mathrm{i})}$ \cite{suppmat}}.
\pc{Combined with Eq. \eqref{coupling matrix D}, this formula implies that}
\begin{eqnarray}
\frac{\gamma_\mathrm{i}}{\gamma}=1-\dfrac{\sigma+\sqrt{\epsilon^2+1}}{2}. \label{growth rate sigma relation}
\end{eqnarray}
\pc{The difference between $\gamma$ and $\gamma_\mathrm{i}$ is the reversible decay rate $\gamma_\mathrm{r}=\gamma-\gamma_\mathrm{i}$. According to Eq. \eqref{growth rate sigma relation}, unitary conditions ($\sigma=1$, $\epsilon=0$) represent a reversible background process ($\gamma=\gamma_\mathrm{r}$). From a physical perspective, for the aeroacoustic scatterer modeled here, $\gamma_\mathrm{r}$ represents radiation losses through the ports and $\gamma_\mathrm{i}$ visco-thermal losses of the pure acoustic system in the absence of a self-oscillating mode.}

\pc{The spectral low-rank approximation used above is exact when $\sigma=1$ and $\epsilon=0$, in which case the spectrum of $\sigma S_*-C$ is singular, and its quality decreases as $|\mu_1|$ and $|\mu_2|$ approach each other \footnote{\pc{In the limiting case $\sigma=\epsilon=0$, $\sigma S_*-C=-C$ has no distinguished dominant spectral subspace.}}.}

\begin{figure}[t!]
\begin{psfrags}
    \psfrag{a}{\hspace{-0.65cm}$|S_{ij}|$ (dB)}
\psfrag{b}{\hspace{-0.05cm}$-29$}
\psfrag{c}{$11$}
\psfrag{d}{\hspace{-0.4cm}$\tilde{s}=9$ }
\psfrag{e}{\hspace{-0.3cm}$|S_{11}|$}
\psfrag{f}{\hspace{-0.3cm}$|S_{12}|$}
\psfrag{g}{}
\psfrag{h}{}
\psfrag{i}{\hspace{1.55cm}\textcolor{red2}{$\leq 1$}}
\psfrag{j}{}
\psfrag{k}{}
\psfrag{l}{}
\psfrag{m}{\hspace{-0.4cm}$\tilde{s}=1.5$ }
\psfrag{n}{\hspace{-0.3cm}$|S_{21}|$}
\psfrag{o}{\hspace{-0.3cm}$|S_{22}|$}
\psfrag{p}{}
\psfrag{q}{}
\psfrag{r}{}
\psfrag{s}{}
\psfrag{t}{}
\psfrag{u}{}
\psfrag{v}{\hspace{-0.4cm}$\tilde{s}=0.6$ }
\psfrag{w}{\hspace{-0.12cm}$\alpha_1$}
\psfrag{x}{\hspace{-0.12cm}$\alpha_2$}
\psfrag{y}{\hspace{-0.375cm}$3.5$}
\psfrag{z}{\hspace{-0.25cm}$10$}
\psfrag{A}{\hspace{-0.45cm}$|S_{ij}|$ (-)}
\psfrag{B}{\hspace{-0.25cm}$\tilde{s}$ (-)}
\psfrag{C}{\hspace{-0.1cm}$0$}
\psfrag{D}{\hspace{-0.3cm}$0.6$}
\psfrag{E}{ \hspace{-0.33cm}\fc{$1.74$}\hspace{1.3cm}\fc{$1.9$}}
\psfrag{F}{}
\psfrag{G}{}
\psfrag{H}{\hspace{-0.2cm}$25$}
\psfrag{I}{}
\psfrag{J}{\hspace{-0.2cm}$25$}
\psfrag{K}{\hspace{-1cm}\hspace{0.5cm}$\dfrac{\omega}{2\pi}$ (arb. units)}
\psfrag{L}{}
\psfrag{M}{$-14$}
\psfrag{N}{$0.5$}
\psfrag{O}{\hspace{-0.3cm}$\alpha_j$ (-)}
\psfrag{P}{\hspace{-0.27cm}$|S_{ij}|$}
\psfrag{1}{{(a)}}
\psfrag{2}{{(d)}}
\psfrag{3}{{(g)}}
\psfrag{4}{{(b)}}
\psfrag{5}{(e)}
\psfrag{6}{{(h)}}
\psfrag{7}{{(c)}}
\psfrag{8}{(f)}
\psfrag{9}{\hspace{0.1cm}(j)}
\psfrag{Z}{}
\psfrag{Y}{\hspace{-0.1cm}\textcolor{red}{L}}
\psfrag{<}{\hspace{-2.5cm}\textcolor{red}{S}\hspace{2.28cm}\textcolor{red}{S}}

    \includegraphics[width=0.45\textwidth]{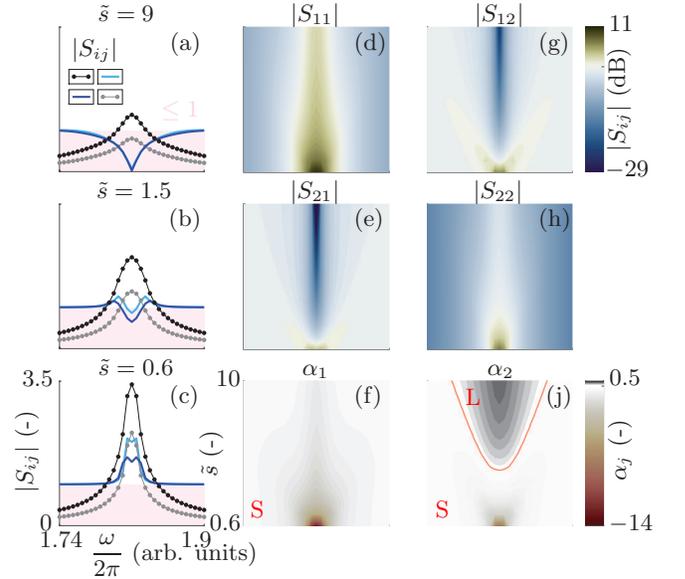}
    \end{psfrags}
    \vspace{-0.2cm}
    \caption{Numerical example of \pc{1D} scattering for \fc{$\omega_0/2\pi=1820$}, $\nu$ equal to $0.4\%$ of $\omega_0$, \pc{$\gamma=2\nu$}, $\sigma=0.6$, $\epsilon=0.3$ and $\kappa=1$. Shown are (a)-(e), (g)\rr{,} (h) the \sr{scattering matrix} coefficients $|S_{ij}|$ and (f), (j) the absorption coefficients $\alpha_j$. The interval $[0, 1]$ in (a)-(c) is colored in red. The red curves in (f), (j) mark the \rr{$\alpha_j=0$} contours\sr{, separating the domains of superradiant (S) and lossy (L) scattering.} The scattering matrix $S$ is evaluated using Eq. \eqref{analytical scattering matrix} for forcing at frequency \pc{$\omega/2\pi$} and normalized amplitude $\tilde{s}=s/\sqrt{\gamma}|a_0|$. \pc{The diagonal and off-diagonal entries of $S$ are reflection and transmission coefficients, respectively.}}
    \label{Figure 2}
\end{figure}

For a general model, the nonlinear scattering matrix can be obtained explicitly from the definition $|s_\mathrm{out}\rangle=S|s_\mathrm{in}\rangle$ and the Moore--Penrose pseudoinverse:
\begin{eqnarray}
S=\dfrac{|s_\mathrm{out}\rangle \langle s_\mathrm{in} |}{\langle s_\mathrm{in}| s_\mathrm{in}\rangle}. \label{general scattering matrix}
\end{eqnarray}
\pc{For simplicity,} we focus on synchronized conditions \pc{achievable} for small enough detuning and large enough forcing amplitude \cite{balanov2009simple,suppmat}. For harmonic excitation with \pc{angular} frequency $\omega$, then, we seek the forced response $\tilde{a}$ as the steady-amplitude solution of Eq. \eqref{modal dynamics} \pc{oscillating at the same frequency}. \sr{By \rr{separately applying} forcing from each port with $|s_\mathrm{in}\rangle=s|j\rangle e^{i\omega t}$, using Eq. \eqref{general scattering matrix}, the} NLC scattering matrix \sr{is found to be}
 \begin{eqnarray}
    S=\underbrace{C+D F^{-1}D^\dag}_{\text{linear}}+\underbrace{\frac{1}{s}D\sum^{\pc{2}}_{j=1}\rho_j e^{i \varphi_j} \langle j|}_{\text{nonlinear}}, \label{Definition scattering matrix}
\end{eqnarray}
\sr{where $\tilde{a}_j=\rho_j e^{i(\omega t + \varphi_j)}$ is the forced response resulting from forcing the cavity through the $j$th port \pc{and $| j\rangle$ is the unit vector in $j$-direction}.} \pc{The frequency-dependent quantities} $\rho_j\in\tr{\mathbb{R}^{+}}$ and $\varphi_j\in \tr{[0,2\pi)}$ satisfy Eq. \eqref{modal dynamics}\tr{, which is equivalent to}
\begin{eqnarray}
0&=&\kappa^2 \rho_j^6-2 \nu \kappa \rho_j^4 +(\nu^2+\Delta^2) \rho_j^2-|D_{j}|^2 s^2, \label{amp. Eq.} \\
\varphi_j&=&-\arg{D_{j}}-\arcsin{\bigg(\dfrac{\Delta \rho_j}{|D_{j}|s}\bigg)},
\label{phase Eq.} \end{eqnarray}
 where \rr{$\Delta=\omega-\omega_0$} is the detuning \pc{and the linearly unstable branch of $\varphi_j$ has been discarded \cite{suppmat}.} Equation \eqref{amp. Eq.} is a cubic equation for $\rho_j^2$ which, for the parameter values considered \tr{in this work}, \sr{always} has a single real \pc{root} \pc{\cite{suppmat,hakmem}.} The analytical $S$-matrix for a biased waveguide \pc{with two ports} coupled to a \sr{nonlinear self-oscillating mode} with $\nu>0$ \sr{follows directly from Eqs. \eqref{coupling matrix D} and \eqref{Definition scattering matrix}}:
\begin{eqnarray}
S&=& \begin{pmatrix*}[c]
\dfrac{\pc{\gamma} g(\epsilon)^2h(\sigma,\epsilon)}{\pc{\gamma} \pc{+} i\Delta} &1- \dfrac{\pc{\gamma} g(\epsilon)h(\sigma,\epsilon)}{\pc{\gamma} \pc{+} i\Delta} \\
1- \dfrac{\pc{\gamma} g(\epsilon)h(\sigma,\epsilon)}{\pc{\gamma} \pc{+} i\Delta}   &\dfrac{\pc{\gamma} h(\sigma,\epsilon)}{\pc{\gamma} \pc{+} i\Delta} 
\end{pmatrix*} \nonumber\\
&&+\dfrac{\sqrt{\pc{\gamma} h(\sigma,\epsilon)}  }{s}\begin{pmatrix*}[c]
-{\rho_1 e^{i\varphi_1} g(\epsilon)}& -{\rho_2 e^{i\varphi_2} g(\epsilon)}\\
{\rho_1 e^{i\varphi_1}} &{\rho_2 e^{i\varphi_2}}
\end{pmatrix*}. \label{analytical scattering matrix}
\end{eqnarray}
At each \pc{frequency} $\omega$, \pc{the forced response is} obtained by \pc{first solving} Eq. \eqref{amp. Eq.} \pc{for $\rho_{1,2}$} and substituting the results into Eq. \eqref{phase Eq.} \pc{to compute $\varphi_{1,2}$. By prescribing different values for $\sigma$ and $\epsilon$ in $S_*$, one can tailor $D$ such that the nonlinear wave-mode coupling leads to superradiant scattering at a given forcing amplitude $s$.}

We now focus on the distribution of energy of the scattered waves in the ports, which is determined by the absolute values of the entries of $S$, $|S_{ij}|$  \fc{\footnote{\fc{Referring to ``energy'' in a generic sense as the squared Euclidian norm; for a detailed discussion of acoustic energy and scattering in the presence of mean flow, see \cite{auregan1999determination}.}}}. We define the absorption coefficients for forcing from the $j$th port as $\alpha_{j}=1-|S_{1 j}|^2-|S_{2 j}|^2$, $j=1,2$. Superradiant scattering is defined as $\alpha_j<0$, implying a net amplification of the incident wave energy\pc{: $\langle s_\mathrm{out}| s_\mathrm{out}\rangle>\langle s_\mathrm{in}| s_\mathrm{in}\rangle$.}

In the \pc{linear,} \fc{symmetric and unitary} limit, the energy reflection coefficient $|S_{11}|^2$ corresponds to the TCMT result of Fan et al., given by Eq. (17) of Ref. \cite{Fan2003569}. For comparison, set $\{ \rho_1,\rho_2,\epsilon,\sigma\}=\{0,0,0,1\}$ in Eq. \eqref{analytical scattering matrix} and $\{1/\tau,\phi,t,r \}=\{\pc{\gamma},-\pi/2,1,0 \}$ in the reference.

A numerical example of the scattering matrix given by Eq. \eqref{analytical scattering matrix} is shown in Fig. \ref{Figure 2} \footnote{\pc{This work uses the scientific color maps of Fabio Crameri described in Ref. \cite{Crameri2020}.}}. $S$ strongly depends on the normalized forcing amplitude $\tilde{s}=s /\sqrt{\gamma} |a_0|$, undergoing a \sr{nonlinear transition from \tr{omnidirectional to purely reflectional superradiant scattering}} as $\tilde{s}$ is increased. The absorption coefficients are each negative over a continuous range of amplitudes and frequencies, resulting from the energy production of the self-sustained mode. While the transition at low amplitudes is sharp, at larger amplitudes\pc{, a} saturation occurs and the superradiant state persists with little change.

To confirm our theoretical analysis, a superradiant scatterer \pc{was} experimentally realized by \tr{means of a self-sustained aeroacoustic mode in a side cavity which is connected to an acoustic waveguide with anechoic terminations. The cavity whistling \pc{was} obtained by imposing a low Mach air flow in the waveguide with a bulk velocity that exceeds the threshold of a supercritical Hopf bifurcation. The corresponding aeroacoustic limit cycle at $\omega_0/2\pi=1.82$ kHz involves a longitudinal mode of the cavity which constructively interacts with the coherent vorticity fluctuations in the cavity opening} \cite{Elder1982532,Howe1980407,Boujo2020,Bourquard2020,Pedergnana2021}. As a byproduct of the air stream, a bias on the order of the Mach number $\mathrm{Ma}\approx 0.17$ is imposed on the system. The measured root-mean-square (rms) \pc{pressure fluctuations} of the self-sustained acoustic waves radiated \pc{into} the upstream and downstream sections of the waveguide are $0.54$ and $0.65$ hPa. The \tr{acoustic} energy of the limit cycle feeds the scattered waves, \tr{enabling} superradiance \tr{in the presence of incident waves}. 

\begin{figure}[t!]
\begin{psfrags}
    
\psfrag{a}{\hspace{-0.65cm}$|S_{ij}|$ (dB)}
\psfrag{b}{\hspace{-0.05cm}$-29$}
\psfrag{c}{\hspace{-0.05cm}$11$}
\psfrag{d}{\hspace{-0.78cm}$s=1.2$ hPa}
\psfrag{e}{\hspace{-0.3cm}$|S_{11}|$}
\psfrag{f}{\hspace{-0.3cm}$|S_{12}|$}
\psfrag{g}{}
\psfrag{h}{}
\psfrag{i}{\hspace{1.55cm}\textcolor{red2}{$\leq 1$}}
\psfrag{j}{}
\psfrag{k}{}
\psfrag{l}{}
\psfrag{m}{\hspace{-0.78cm}$s=0.3$ hPa}
\psfrag{n}{\hspace{-0.3cm}$|S_{21}|$}
\psfrag{o}{\hspace{-0.3cm}$|S_{22}|$}
\psfrag{p}{}
\psfrag{q}{}
\psfrag{r}{{}}
\psfrag{s}{}
\psfrag{t}{}
\psfrag{u}{}
\psfrag{v}{\hspace{-0.78cm}{$s=0.1$} hPa}
\psfrag{w}{\hspace{-0.12cm}$\alpha_1$}
\psfrag{x}{\hspace{-0.12cm}$\alpha_2$}
\psfrag{y}{\hspace{-0.375cm}$3.5$}
\psfrag{z}{\hspace{-0.35cm}$1.3$}
\psfrag{A}{\hspace{-0.45cm}$|S_{ij}|$ (-)}
\psfrag{B}{\hspace{-0.45cm}$s$ (hPa)}
\psfrag{C}{\hspace{-0.1cm}$0$}
\psfrag{D}{\hspace{-0.3cm}$0.1$}
\psfrag{E}{ \hspace{-0.33cm}$1.74$\hspace{1.3cm}$1.9$}
\psfrag{F}{\hspace{-0.3cm}$1.9$}
\psfrag{G}{}
\psfrag{H}{\hspace{-0.35cm}$1.3$}
\psfrag{I}{}
\psfrag{J}{\hspace{-0.2cm}$1.3$}
\psfrag{K}{\hspace{-0.61cm}$\dfrac{\omega}{2\pi}$ (kHz)}
\psfrag{L}{}
\psfrag{M}{$-14$}
\psfrag{N}{$0.5$}
\psfrag{O}{\hspace{-0.3cm}$\alpha_j$ (-)}
\psfrag{P}{\hspace{-0.27cm}$|S_{ij}|$}
\psfrag{1}{(a)}
\psfrag{2}{(d)}
\psfrag{3}{{(g)}}
\psfrag{4}{(b)}
\psfrag{5}{(e)}
\psfrag{6}{(h)}
\psfrag{7}{(c)}
\psfrag{8}{(f)}
\psfrag{9}{\hspace{0.05cm}(j)}
\psfrag{W}{\hspace{0.13cm}\textcolor{red}{S}}
\psfrag{Z}{\hspace{-0.1cm}\textcolor{red}{L}\hspace{2.245cm}\textcolor{red}{L}}
\psfrag{Y}{\hspace{-0.015cm}\textcolor{red}{S}}
\psfrag{<}{\textcolor{red}{S}}

    \includegraphics[width=0.45\textwidth]{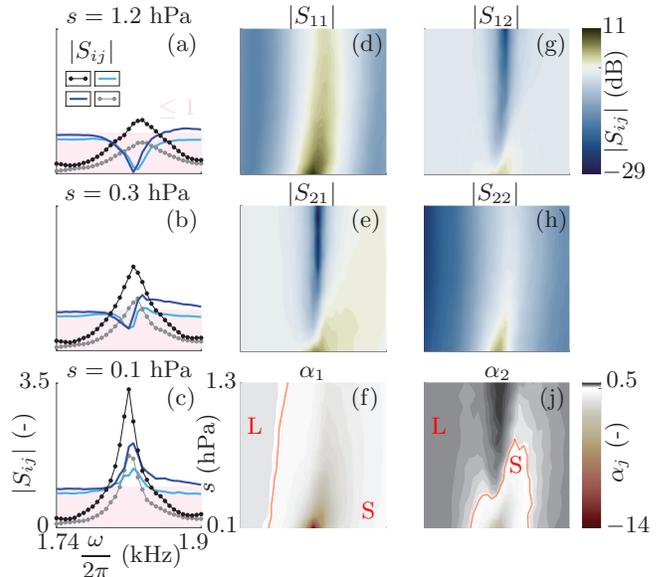}
    \end{psfrags}
        \vspace{-0.2cm}
    \caption{Experimentally determined scattering matrix \tr{of a side cavity in} a biased waveguide \pc{with two ports}\rr{, when the cavity exhibits a self-oscillating aeroacoustic mode} at $\omega_0/2\pi= 1.82$ kHz. Shown are the same quantities as in Fig. \ref{Figure 2}. To measure $S$, acoustic forcing at frequency \pc{$\omega/2\pi$} with amplitude $s$ is applied \tr{from the two anechoic ends of the waveguide.} \tr{The rms \rr{amplitudes} of the acoustic waves radiated by the cavity whistling in absence of acoustic forcing are $0.54$ and $0.65$ hPa in the \pc{upstream} and downstream sections of the waveguide.}}
    \label{Figure 3}
\end{figure}

To measure the scattering matrix $S$, acoustic \tr{waves produced by compression drivers placed in the anechoic ends of the waveguide are sent to the cavity.} The columns $S_{i1}$ and $S_{i2}$, $i=1,2$ are obtained separately by applying \tr{harmonic forcing} \pc{upstream} and downstream of the cavity. The \tr{voltage of the signal} to the drivers is calibrated at each frequency to achieve a specified \tr{acoustic} forcing amplitude $s$ in the waveguide \pc{just outside} the cavity aperture. Assuming lossless transmission, the multi-microphone method is employed to reconstruct the forward- and back\pc{ward}-propagating waves in the \tr{waveguide} \cite{Jang19981520}. Details of the experimental set-up are provided in the Supplemental Material \cite{suppmat}.

\rr{Similar to the numerical example presented \pc{in Fig. \ref{Figure 2}},} the experimentally determined $S$-matrix also \rr{exhibits} a nonlinear \tr{transition} with increasing $s$ from \tr{omnidirectional to purely reflectional superradiant scattering} (Fig. \ref{Figure 3}). It is interesting to note that in the present biased system, superradiant reflection is stronger for incident waves in flow direction and persists at large forcing amplitudes, confirming its theoretically predicted robustness.

\begin{figure}[t!]
\begin{psfrags}
    
\psfrag{w}{\hspace{-0.55cm}NLC off}
\psfrag{v}{\hspace{-0.55cm}NLC on}
\psfrag{y}{\hspace{-0.35cm}$1.5$}
\psfrag{A}{\hspace{-0.4cm}$|S_{ij}|$ (-)}
\psfrag{G}{$1.74$\hspace{1.4cm}$1.9$}
\psfrag{E}{}
\psfrag{K}{}
\psfrag{L}{\hspace{-0.6cm}$\dfrac{\omega}{2\pi}$ (kHz)}
\psfrag{C}{\hspace{-0.04cm}$0$}
\psfrag{P}{\hspace{-0.05cm}$|S_{ij}|$}
\psfrag{Z}{\hspace{0.1cm}\textcolor{red2}{$\leq 1$}}
    \includegraphics[width=0.32\textwidth]{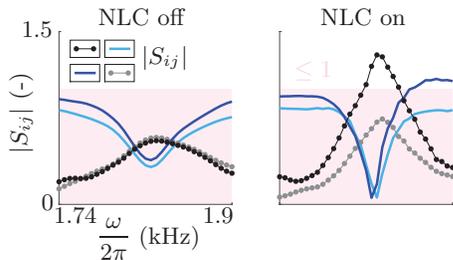}
    \end{psfrags}
            \vspace{-0.3cm}
    \caption{Experimentally determined \sr{scattering matrix} coefficients $|S_{ij}|$ at forcing amplitude $s=1.3$ hPa \pc{for resonance-based (left inset) and limit cycle-amplified scattering (right inset).} The interval $[0, 1]$ is colored in red.}
    \label{Figure 4}
\end{figure}

In the unbiased case, for $\mathrm{Ma}=0$ (without air \tr{flow}), the scattering properties of the resonant cavity are \pc{amplitude-independent}, \pc{quasi-}symmetric and \tr{lossy} (Fig. \ref{Figure 4}\pc{, left inset}). The \pc{scattering matrix} in the unbiased case is presented in detail in the Supplemental Material \cite{suppmat}.

\fc{The self-sustained aeroacoustic system considered in this work can be compared to an electromagnetic circuit with a negative resistance element, such as a Gunn diode \cite{Knight1967393,Mosekilde19902298} or a tunnel diode \cite{Bernard1963627,Nguyen2013,Fluckey2022}, except that the negative resistance in our case does not come from an external element but from the aeroacoustic instability which is intrinsic to the medium that the waves travel in.}

\fc{We investigated both small-amplitude amplification by a Hopf bifurcation \cite{Martin19864523,Lacot200310} (see Figs. \ref{Figure 2}c and \ref{Figure 3}c) as well as scattering in the strongly nonlinear regime (see Figs. \ref{Figure 2}a and \ref{Figure 3}a). Our tunable concept is based on the synchronization between incident waves and the radiation losses of a self-oscillating source. Like in the nonlocal active metamaterial proposed in \cite{Geib2021}, destructive interferences between traveling waves in absence of the source and waves radiated from the source enable tailoring the scattering to quasi-perfect isolation. At the operating condition shown in Fig. \ref{Figure 4} (right inset), the scatterer is a reflectional amplifier for waves originating from port ``1'' \cite{Westig2018,Zhong2020}.} \tr{We \pc{expect} that the quasi-passive realization of superradiance presented in this letter will find application} in \rr{flow-based acoustic \fc{metamaterials} \fc{\cite{Cummer2016,Fleury2014516,Auregan2017,Ding2019}}, which often suffer from dissipation losses}.

\fc{Future applications of our theory may also include modeling of optical scatterers such as complex nonreciprocal devices \cite{Xu2015,Zhong2020} and spherical cavities \cite{Xie2022}. While the derivation of the background could be performed analogously by choosing suitable matrices $S_*$ and $C$, nonlinearities \cite{StephenYeung19984421,Roke2004115106,Shim2016,Silver2022}, negative feedback \cite{DeAssis2017} or thermal noise \cite{Webster2008} may be included through the modal equation(s). The scattering matrix follows directly from the general expression \eqref{general scattering matrix}.}

This project is funded by the Swiss National Science Foundation under Grant agreement 184617. 

\bibliographystyle{apsrev4-2} 
\bibliography{apssamp}
\end{document}